\documentclass[12pt,reqno]{article}

\usepackage{amssymb}
\usepackage{graphicx}
\usepackage{epstopdf}

\usepackage{fullpage}

\usepackage{graphics,amsmath,amssymb}
\usepackage{amsthm}
\usepackage{amsfonts}

\usepackage{enumerate}

\usepackage[hyphens]{url}

\def\keywords{\vspace{.5em}
{\textit{Keywords}:\,\relax%
}}

\begin{document}

\newtheorem{thm}{Theorem}[section]
\newtheorem{proposition}[thm]{Proposition}
\newtheorem{definition}[thm]{Definition}

\theoremstyle{remark}
\newtheorem{remark}[thm]{Remark}

\begin{center}
\vskip 1cm{\Large\bf On property of least common multiple to be a $D$-magic number}
\vskip 1cm
{\small  V.L. Gavrikov\\
Institute of Ecology and Geography \\
Siberian Federal University\\
600041 Krasnoyarsk, pr. Svobodnyi 79\\
Russian Federation\\
vgavrikov@sfu-kras.ru
}
\end{center}

\begin{abstract}
Least common multiple ($lcm$) has been shown to posses the property of $D$-magic number, that is, its least significant digit $0$ does not change when the number is transferred into all other numbering systems with smaller bases. The number $lcm + 1$ preserves this property as well.

\end{abstract}

\section{Introduction}
Least common multiple ($lcm$) is a function which was often referred to as having two arguments, i.e. $lcm[x_1, x_2]$ but can be easily reformulated to any number of arguments, $lcm[x_1, x_2, \dots, x_n]$.

The function has been widely known for being used at formulating of encryption algo\-rithms, both in classical works \cite{shannon1949} and in later research on encryption keys \cite{rivest1978}. Because of its important applications properties of $lcm[\:]$ are of interest. An identity has been proven \cite{farhi2009} that relates $lcm[\:]$ of binomial coefficients to $lcm[\:]$ of the sequence of indices of the coefficients. A typical behavior of $lcm[\:]$ of random subsets \{1, \dots, n\} \cite{cilleruelo2014} has also been studied.

In this work, some properties of divisibility of $lcm[\:]$ function are explored that lead to a sort of invariance of the least significant digit of a number when the number is transferred to a different numbering system.

As usual, when a multidigit integer is transferred to a numbering system its least significant digit (as well as other digits) changes, e.g., $64_{10} = 100_8, 100_{10} = 244_6$. Sometimes however the transfer to another numbering sysytem does not lead to the change in the least significant digit, e.g., $126_{10} = 176_8, 101_{10} = 401_5$.

From these observations let us put a more general question: how can one get the the number that does not change its least significant digit when being transferred to another numbering system?

\section{Formulation}

\begin{definition}
For an arbitrary base-$L$ numbering system, $D$-magic number $M$ is such a numberthat does not change its least significant digit when being transferred to any other base-$l$ numbering system, with $l < L$.
\end{definition}

An integer number $M$ in base-$L$ system may be represented in \textit{decimal} form:
\begin{equation}
M_L = L \cdot n + j,
\label{eq:1}
\end{equation}
where $n$ is the number of tens in $M_L$ and $j$ is the least significant digit of $M_L$, with $j < L$.

If $l$ is the base of numbering system then the transfer from $M_L$ to $M_l$ will include calculations of remainders from division by $l$ both $L \cdot n$ and $j$. Provided these remainders are known a new value for $j$ is received.

If $L \cdot n$ in Eq. (\ref{eq:1}) is divisible without a remainder by all $l$, $2 \leq l < L$, and $j < l$ then $j$ will not change when $M_L$ is transferred to any base-$l$ system. There is an infinite quantity of numbers divisible by all $2 \leq l < L$ but the minimal of them is ony one. And this number is least common multiple. In other words, $lcm[\forall l, 2 \leq l < L]$ is a $D$-magic number in base-$L$ system (as well as in all systems with bases smaller than $L$). Therefore, \textit{calculation of $lcm[\:]$ is the very algorithm to get $D$-magic numbers}.

\section{Illustrations}

It is easy to find, e.g., in base-ten system , such a number that wil be divisible without a remainder by $10, 9, 8, 7, 6, 5, 4, 3, 2$. As well known, $lcm[10, 9, 8, 7, 6, 5, 4, 3, 2]$ = $2520$ (see sequence A003418 in On-line Encyclopedia of Integer Sequences OEIS \sloppy \url{http://oeis.org/A003418}).

\fussy

A transfer of decimal number $2520$ to any numbering system with bases $l < 10$ does not change the least significant digit (in this particular case  $j = 0$):\\
\begin{center}
\begin{tabular}{lr}
$l$&$M_l$\\
    10    & 252\textbf{0}\\
    9     & 341\textbf{0}\\
    8     & 473\textbf{0}\\
    7     & 1023\textbf{0}\\
    6     & 1540\textbf{0}\\
    5     & 4004\textbf{0}\\
    4     & 21312\textbf{0}\\
    3     & 1011010\textbf{0}\\
    2     & 10011101100\textbf{0}\\
\end{tabular}\\
\end{center}

Moreover, in case $j = 1$ (see Eq. \ref{eq:1}) this least significant digit will not change as well:\\
\begin{center}
\begin{tabular}{lr}
$l$&$M_l$\\
    10    & 252\textbf{1}\\
    9     & 341\textbf{1}\\
    8     & 473\textbf{1}\\
    7     & 1023\textbf{1}\\
    6     & 1540\textbf{1}\\
    5     & 4004\textbf{1}\\
    4     & 21312\textbf{1}\\
    3     & 1011010\textbf{1}\\
    2     & 10011101100\textbf{1}\\
\end{tabular}\\
\end{center}

If $j \in \{2, 3, 4, 5, 6, 7, 8\}$ such a property (constance of least significant digit) holds only at $j < l$.

\begin{remark}
Thus $lcm[10, 9, 8, 7, 6, 5, 4, 3, 2]$ equal to $2520$ not only is $D$-magic number itself for base-ten numbering system but also produces a set of $D$-magic numbers--by adding of least significant digit $j < 10$.
\end{remark}

Let us now look at how this approach works at $L \not= 10$.

For base-eight system, $lcm[8, 7, 6, 5, 4, 3, 2]$ = $840_{10}$ = $1510_8$. It can be seen that base-eight number $1510$ does not change least significant digit when being transferred into numbering systems with bases $7, 6, 5, 4, 3, 2$:\\
\begin{center}
\begin{tabular}{lr}
$l$&$M_l$\\
    8     & 151\textbf{0}\\
    7     & 231\textbf{0}\\
    6     & 352\textbf{0}\\
    5     & 1133\textbf{0}\\
    4     & 3102\textbf{0}\\
    3     & 101101\textbf{0}\\
    2     & 110100100\textbf{0}\\
\end{tabular}\\
\end{center}

Correspondingly, the base-eight number $1511_8$ will also not change least significant digit when transferred into system with bases smaller than $8$.

Another example, base-16 numbering system. $lcm[16, 15, 14, 13, 12, 11, 10, 9, 8, 7, 6, 5, 4,$ $ 3, 2]$ = $720720_{10}$ = $aff50_{16}$. The transfer of number $aff50_{16}$ into systems with bases smaller than $16$ gives:\\
\begin{center}
\begin{tabular}{lr}
$l$&$M_l$\\
   16     & aff5\textbf{0}\\
   15     & e383\textbf{0}\\
   14     & 14a92\textbf{0}\\
   13     & 1c308\textbf{0}\\
   12     & 2a910\textbf{0}\\
   11     & 45254\textbf{0}\\
   10     & 72072\textbf{0}\\
    9     & 131757\textbf{0}\\
    8     & 257752\textbf{0}\\
    7     & 606114\textbf{0}\\
    6     & 2324040\textbf{0}\\
    5     & 14103034\textbf{0}\\
    4     & 223333110\textbf{0}\\
    3     & 110012112210\textbf{0}\\
    2     & 1010111111110101000\textbf{0}\\
\end{tabular}\\
\end{center}

Therefore least common multiple of $2, 3,$ $ 4, 5 \dots L$ is a $D$-magic number for the numbering system with the base $L$ (maximum of this sequence). A convenient algorithm could be as follows: 1) first, one gets $lcm[2, 3, 4, 5 \dots L]$ for base-ten system and then 2) transfers it into system with the base $L$. This procedure leads to the number having $0$ as least significant digit.

Adding of unity $1$ to the list significant digit $0$ brings about another $D$-magic number. Adding of a digit $j \in \{2, 3, 4, 5 \dots L - 1 \}$ to the least significant digit produces a set of set number that are partly $D$-magic; when being transferred into base-$l$ systems the least significant digit $j$ of them will not change only when  $j < l$.

\end{document}